

\baselineskip=14pt
\parskip=10pt

\magnification=\magstephalf

\def\1{{\overline{1}}}
\def\2{{\overline{2}}}
\parindent=0pt
\overfullrule=0in

\def\frac#1#2{{#1 \over #2}}
\centerline
{\bf IDENTITIES In CHARACTER TABLES Of ${\bf S_n}$ }
\bigskip
\centerline
{\it Alon REGEV, Amitai REGEV,  and Doron ZEILBERGER}

{\bf Abstract}: 
{\it We use the algebra of {\it difference operators} to study sums of squares (and other powers) of the
characters of the symmetric group, $\chi^{\lambda}(\mu)$, when the sum is restricted over shapes, $\lambda$, with a fixed number of rows, and for
hook shapes, and $\mu$ has `mostly' ones.
We prove that such sums are always P-recursive, i.e., satisfy a linear difference equation with polynomial coefficients.
For the special case of two rows, and for hook-shapes, we prove that these sums are in fact closed-form, and we present
algorithms, complete with rigorous proofs, for finding these expressions.This article is accompanied by a Maple
package, {\tt Sn} (available from {\tt http://www.math.rutgers.edu/\~{}zeilberg/tokhniot/Sn}), and a webpage,
{\tt http://www.math.rutgers.edu/\~{}zeilberg/mamarim/mamarimhtml/sn.html}, with links to extensive output,
containing rigorously-proved explicit formulas for many cases.}

{\bf Introduction}

The notion of {\it group character}, introduced by Frobenius,
is a fundamental one in representation theory, and among all groups, the symmetric group, $S_n$, is the
most fundamental one. The entries of the character tables of the symmetric group
may be defined in numerous ways, but for our purposes, we will use
the following definition, that does not require any knowledge of `advanced' algebra, high-school algebra suffices!

Recall that the {\it Constant Term} of a {\it Laurent polynomial} in $(x_1, \dots, x_m)$ is the free term, i.e.
the coefficient of $x_1^0 \cdots x_m^0$. For example

$$
CT_{x_1,x_2} (x_1^{-3}x_2+x_1x_2^{-2}+5)=  5 \quad  .
$$

Recall that a {\it partition} (alias {\it shape}) of an integer $n$, with $m$ {\it parts} (alias {\it rows}),
is a non-increasing sequence of positive integers
$$
\lambda= (\lambda_1, \dots, \lambda_m) \quad,
$$
where $\lambda_1 \geq \lambda_2 \geq \dots \geq \lambda_m >0$, and $\lambda_1 + \dots + \lambda_m=n$.

If $\lambda=(\lambda_1, \dots, \lambda_m)$ and $\mu=(\mu_1, \dots, \mu_r)$ are partitions of $n$ with $m$ and $r$ parts,
respectively, then it easily follows from (7.8) (p. 114) in [M], that
the {\it characters}, $\chi^{\lambda}(\mu)$,  of the {\it symmetric group}, $S_n$, may be obtained
by the {\it constant term} expression
$$
\chi^{\lambda}(\mu) \, = \,
CT_{x_1, \dots , x_m} \,
\frac
{\prod_{1 \leq i < j \leq m} (1-\frac{x_j}{x_i}) \prod_{j=1}^r \left ( \sum_{i=1}^{m} x_i^{\mu_j} 
\right )}
{\prod_{i=1}^{m} x_i^{\lambda_i}} 
  \quad .
\eqno(Chi)
$$

It is well-known (e.g. [M], p. 119, bottom line) that, writing $\mu$ in {\it frequency notation}, $\mu=1^{a_1} 2^{a_2} \dots n^{a_n}$,
we have the following beautiful identity:
$$
\sum_{\lambda \vdash n}  \chi^{\lambda}(\mu)^2 \, = \, \prod_{i=1}^{n} i^{a_i} a_i! \quad.
$$
The most famous special case is when $\mu=1^n$, that becomes 
(since $\chi^{\lambda}(1^n)=f_\lambda$, the number of Standard Young Tableaux of shape $\lambda$),
the identity:
$$
\sum_{\lambda \vdash n}  f_{\lambda}^2 \, = \, n! \quad,
$$
that has a lovely combinatorial proof using the celebrated Robinson-Schensted correspondence ([Ro][Sc]).
A bijective proof of the former, more general, identity was given by Dennis White ([Wh]).

If one restricts the latter sum to go over partitions with at most a fixed number of parts, then
one gets many sequences of combinatorial interest. Most notably, for a fixed $r \geq 1$,
$$
a^{(r)}(n)\, := \, \sum_{ { {\lambda \vdash n} \atop {length(\lambda) \leq r}}}  f_{\lambda}^2 \quad,
$$
is the number of permutations of length $n$ that do not contain an increasing subsequence of length $r+1$.
The sequences $a^{(r)}(n)$ are (as of July 8, 2015) in [Sl] for $r \leq 11$. Note that $a^{(2)}(n)$ is the super-famous sequence 
{\bf A000108}, of Catalan numbers, $C_n:=\frac{1}{n+1} {{2n} \choose {n}}$.

In [Z1] it was proved that, for any fixed $r$, the sequence $a^{(r)}(n)$ is $P$-recursive (alias {\it holonomic}), i.e. satisfies
a {\it linear difference equation} (alias {\it linear recurrence equation}) with {\it polynomial coefficients}.
We should also mention that Ira Gessel ([Ge]) famously discovered a lovely determinant formula for the generating functions in terms
of Bessel functions.

The analogous sequences, for the straight sums (without the squares),
$$
b^{(r)}(n) \, := \,\sum_{ { {\lambda \vdash n} \atop {length(\lambda) \leq r}}}  f_{\lambda} \quad,
$$
are also of combinatorial interest, counting the number of {\it involutions} of length $n$ avoiding
an increasing subsequence of length $r+1$. It is well known,  and easy to prove, that $b^{(2)}(n)={{n} \choose {\lfloor n/2 \rfloor}}$,
that is sequence {\bf A001405} in [Sl]. Much deeper is the result, first proved in [Re],
that $b^{(3)}(n)$ are the {\it Motzkin numbers}, {\bf A001006} (see [Z2] for another proof and for
a generalization). The sequence $b^{(4)}(n)$, {\bf A005817},
is even nicer, given in terms of Catalan numbers, $b^{(4)}(n)=C_{\lfloor n/2 +1/2\rfloor} C_{\lfloor n/2 +1 \rfloor}$,
as first proved by Dominique Gouyou-Beauchamps ([Go]). $b^{(5)}(n)$ is {\bf A049401}, while
$b^{(6)}(n)$ is {\bf A007579}. See also [BFK].

Recall ([Z1]) that a discrete function $a(n)$ of a {\it single} variable is called {\it holonomic} (or $P$-recursive)
if it satisfies a homogeneous linear recurrence (alias {\it difference}) equation with polynomial coefficients, i.e.
there exists an integer $L$ and polynomials $p_0(n), \dots , p_L(n)$, such that
$$
\sum_{i=0}^{L} p_i(n) a(n+i) \, = 0 \quad, \quad (n \geq 0) \quad .
$$
A discrete function of several variables $a(n_1, \dots, n_m)$ is holonomic if it satisfies such a recurrence
in each of its  variables, $n_1, \dots, n_m$ and the coefficients are polynomials in all of them, and the system is
non-degenerate. It was proved in [Z1] that if you sum such a holonomic discrete function over some of its arguments,
you get yet-another holonomic function in the surviving variables.

As mentioned in [Z1], the reason that $a^{(r)}(n)$ and $b^{(r)}(n)$ (and more generally, sums of powers $f_{\lambda}^s$, 
where $s$ is any positive integer), for any {\it specific}, $r$, are
$P$-recursive in $n$ is that the {\it summand}, $f_{\lambda}=f(\lambda_1, \dots, \lambda_r)$, and hence any of its powers,
is {\it holonomic} in its arguments, and hence the $(r-1)$-fold {\it multisum} is guaranteed to be holonomic in
the `surviving' discrete variable $n$. Furthermore, thanks to the {\it hook-length} formula ([Wi]), or equivalently
the {\it Young-Frobenius} formula, the summand $f(\lambda_1, \dots, \lambda_r)$  can be expressed as
$$
f(\lambda_1, \dots, \lambda_r) \, = \,
\frac{\prod_{1 \leq i <j \leq r} (\lambda_ j-\lambda_i+j-i) (\sum_{i=1}^{r} \lambda_i)!}
{\prod_{i=1}^{r} (\lambda_i+r-i)!} \quad.
$$
In addition, thanks to [AZ], there are effective algorithms for finding these recurrences, implemented in the
Maple package  \hfill\break
{\tt http://www.math.rutgers.edu/\~{}zeilberg/tokhniot/MultiZeilberger} \quad .

However, it is much more efficient to derive these recurrences by generating sufficiently many terms, and
then {\it guessing} the linear recurrence, that we know for sure exists, by {\it undetermined coefficients}.
These recurrences can be proved fully rigorously, if desired, but since they are definitely true, we do not
bother to waste time on this.

More generally, one can consider such sums where the shape $\lambda$ belongs to a ``meta-hook'' with $k$ rows and
$l$ columns, in other words, the analogous sums (see [BR] and [EZ]) where  one sums $f_\lambda$ (or its square,
or any positive integer power), over all shapes, that do {\bf not} contain the celll $(k+1,l+1)$. Such a shape
is determined by $k+l$ parameters, the lengths of the $k$ largest rows and the lengths of the
$l$ largest columns, and, once again, using the Hook Length Formula, one can express $f_\lambda$ as
a closed-form expression in these $k+l$ discrete parameters, and the above observations about
the resulting sequences, for each {\it fixed} $(k,l)$, being holonomic, still apply.

{\bf Character Sums for $\mu$ `Close' to $1^n$}

For any partition $\lambda$, let $|\lambda|$ be its sum, in other words, the integer that is being partitioned.

As we noted above, $f_\lambda$ equals  $\chi^{\lambda}(1^n)$. The first  purpose of the
present paper is merely to {\bf observe}, that an analogous argument still applies
if one replaces $\mu=1^n$ by $\mu=\mu_0 1^{n-|\mu_0|}$, for {\it any} fixed partition $\mu_0$ (with smallest part at least $2$),
the analogous sums with $f_\lambda$ replaced by $\chi^{\lambda}(\mu_0 1^{n-|\mu_0|})$ are also {\it guaranteed} to
be holonomic. This follows from Eq. $(Chi)$, that spells out to be, writing $\mu_0=(a_1, \dots, a_k)$ ($a_k \geq 2$).

$$
\chi^{\lambda}(\mu_0 1^{n-|\mu_0|}) \, = \,
CT_{x_1, \dots , x_m} \,
\frac
{\prod_{1 \leq i < j \leq m} (1-\frac{x_j}{x_i}) \prod_{j=1}^k \left ( \sum_{i=1}^{m} x_i^{a_j} \right )
\cdot \left (  \sum_{i=1}^{m} x_i \right )^{n-a_1- \dots -a_k}}
{\prod_{i=1}^{m} x_i^{\lambda_i}}   \quad .
\eqno(Chi0)
$$
Expanding
$$
\prod_{1 \leq i < j \leq m} (1-\frac{x_j}{x_i})
\prod_{j=1}^k \left ( \sum_{i=1}^{m} x_i^{a_j} \right ) \quad,
$$
we get, for {\it fixed} $m$, $k$, and $\mu_0=(a_1, \dots, a_k)$,  a {\bf finite} 
sum of monomials, and collecting the contributions of each, we get
a {\it finite} linear combination of shifts of the  multinomial coefficient $(\lambda_1 + \dots + \lambda_m)!/(\lambda_1! \cdots \lambda_m!)$,
and it is easy to see that the result is a rational function times the latter (with `nice' denominator), and hence closed-form.

The analogous argument for sums over shapes contained in a meta-hook $H(k,l)$ is slightly more complicated, and is omitted.

The second purpose of the present paper is to observe that,  for the special cases
of shapes with at most $2$ rows, defining
$$
\psi^{(2)}(\mu) \, := \, \sum_{j=0}^{\lfloor n/2 \rfloor} \chi^{(n-j,j)}(\mu)^2 \quad,
$$
for {\it any} fixed partition $\mu_0$ (with smallest part $\geq 2$), there is a closed-form expression for
$\psi^{(2)}(\mu_0 1^{n-|\mu_0|})$, of the form
$$
\psi^{(2)}(\mu_0 1^{n-|\mu_0|})\, =\, R_{\mu_0}(n) {{2n} \choose {n}} \quad ,
$$
for some {\bf rational function} $R_{\mu_0}(n)$. The reason is that in the two-rowed case,
$ \chi^{(n-j,j)}(\mu_0 1^{n-\mu0})$ can be expressed (thanks to Eq. $(Chi)$) as a linear combination
of {\it shifts} of the binomial coefficients ${{n} \choose {j}}$, and squaring it and expanding,
gives (possibly many, but still finitely-many) expressions of the form
$$
\sum_{j=0}^{\lfloor n/2 \rfloor}  {{n-\alpha} \choose {j-\beta}} {{n-\alpha'} \choose {j-\beta'}} \quad ,
\quad (\alpha+\alpha'=|\mu_0|) \quad ,
$$
each of which (after symmetrizing in order to make the summations range over all $-\infty<j<\infty$ [with the usual
convention that ${{a} \choose {b}}$ is $0$ if $a<b$ and if $b<0$])
is summable by the Vandermonde-Chu convolution 
$$
\sum_{j}  {{n-\alpha} \choose {j-\beta}} {{n-\alpha'} \choose {j-\beta'}} =  {{2n-\alpha -\alpha'} \choose {n+\beta'-\alpha'-\beta}} , \quad
$$
each of which is a multiple of ${{2n} \choose {n}}$ by a certain rational function, and adding these finitely (but possibly
numerous) terms, still adds up to a certain rational function times   ${{2n} \choose {n}}$.

Analogously, for shapes inside the $(1,1)$-meta hook
$$
\phi^{(2)}(\mu) \, := \, \sum_{j=1}^{n} \chi^{(j,1^{n-j})}(\mu)^2 \quad,
$$
since $\chi^{(j,1^{n-j})}(1^n)={{n-1} \choose {j-1}}$ (as follows easily from $(Chi)$ specialized to this case),
$\chi^{(j,1^{n-j})}(\mu_0 1^{n-|\mu_0|})$ can be expressed as a finite linear combinations of
${{n-1-\alpha} \choose {j-1-\beta}}$, once again we get a finite linear combination of Vandermonde-Chu convolutions, each of
them being a multiple of ${{2n-2} \choose {n-1}}$ by a certain rational function
(equivalently, we can use ${{2n} \choose {n}}$, as above, but it is more natural to use the former as the ``base'', since
it is the answer for $\mu=1^n$, i.e. where $\mu_0$ is the {\it empty} partition).

In fact, in this case we can get nice explicit expressions for the generation functions
(recall that $\mu_0=(a_1, \dots, a_r)$)
$$
F_{n;(a_1, \dots, a_r)}(x):=
\sum_{j=1}^{n} \chi^{(j,1^{n-j})} (\mu_0 1^{n-|\mu_0|}) \, x^j \,=\,
x (1+x)^{n-1-|\mu_0|} \prod_{i=1}^{r} (x^{a_i}-(-1)^{a_i}) \quad .
$$
Hence 
$$
\phi^{(2)}(\mu_0 1^{n-|\mu_0|}) \, = \, \sum_{j=1}^{n} \chi^{(j,1^{n-j})}(\mu_0 1^{n-|\mu_0|} )^2 \quad
$$
is the {\it constant term} of
$$
F_{n;(a_1, \dots, a_r)}(x) \cdot  F_{n;(a_1, \dots, a_r)}(x^{-1}) =
$$
$$
x (1+x)^{n-1-a_1- \dots -a_r} 
\left ( \, \prod_{i=1}^{r} ( x^{a_i}-(-1)^{a_i})  \, \right )
\cdot
x^{-1} (1+x^{-1})^{n-1-a_1- \dots -a_r} 
\left ( \, \prod_{i=1}^{r} ( x^{-a_i}-(-1)^{a_i} )  \, \right )
\quad .
$$
$$
=\frac{(1+x)^{2n-2 -2|\mu_0|}}{x^{n-1-|\mu_0|}} \cdot Q(x) \quad,
$$
where $Q(x)$ is the symmetric Laurent polynomial
$$
Q(x)=\prod_{i=1}^{r} (x^{a_i} -(-1)^{a_i}) (x^{-a_i} - (-1)^{a_i} ) \quad .
$$
Expanding $Q(x)$ as a sum of monomials and extracting the respective coefficients, we get a linear combination
of terms of the form ${{2(n-1-|\mu_0|)} \choose {n-1-|\mu_0|-j}}$, that obviously simplifies to
a rational function times ${{2n-2} \choose {n-1}}$.
[This is implemented in procedure {\tt Phi2} in the Maple package {\tt Sn}].

{\bf Implementation}

Everything discussed here is implemented in the Maple package {\tt Sn}, available from the front of the
present paper

{\tt http://www.math.rutgers.edu/\~{}zeilberg/mamarim/mamarimhtml/sn.html} \quad ,

where one can find several sample input and output files, that readers are welcome to extend.

{\bf Some Output}

I. {\bf Explicit Expressions for} ${\bf \phi^{(2)}_n (\mu_0 1^{n-|\mu_0|})}$

The classical case is well-known, and easy.

$$
\phi^{(2)}_n (1^n)= {{2n-2} \choose {n-1}} \quad .
$$

We also have
$$
\phi^{(2)}_n (2 1^{n-2})= \frac{1}{2n-3} {{2n-2} \choose {n-1}} \quad (=2C_{n-2})  \quad ,
$$
$$
\phi^{(2)}_n (2 2 1^{n-4})=  \frac {3}{ \left( 2\,n -3 \right)  \left( 2\,n-5 \right)}  {{2n-2} \choose {n-1}} 
\quad .
$$
More generally, for $r\geq 0$, we have:
$$
\phi^{(2)}_n (2^r 1^{n-2r})=  \frac{(2r)!(2n-2r-2)!}{r!(n-1)!(n-r-1)!} 
\quad .
$$
Also
$$
\phi^{(2)}_n (3 1^{n-3})= \frac{n^2 \, - \, 7\,n \, + \, 18}{ 4 \left(2\,n -3 \right)  \left( 2\,n-5 \right) }
 {{2n-2} \choose {n-1}}  \quad ,
$$
$$
\phi^{(2)}_n (4 1^{n-4})= \frac {n^2-9\,n+23}{  \left( 2\,n -3\right) \left( 2\,n-5 \right)  \left( 2\,n-7 \right)  }
 {{2n-2} \choose {n-1}}  \quad ,
$$
$$
\phi^{(2)}_n (5 1^{n-5})= 
\frac {{n}^{4} -22\,{n}^{3} +239\,{n}^{2} -1298n +2760} {  16 \left( 2\,n -3 \right) \left( 2\,n-5 \right)  \left( 2\,n-7 \right)  \left( 2\,n-9 \right)  }
 {{2n-2} \choose {n-1}}  \quad ,
$$
$$
\phi^{(2)}_n (32 1^{n-5})= 
\frac {{n}^{2}-15\,n + 74}{ 4 \left( 2\,n -3 \right) \left( 2\,n-5 \right)  \left( 2\,n-7 \right)  }
 {{2n-2} \choose {n-1}}  \quad .
$$

For {\it all} the (proved!) explicit expressions for $\phi^{(2)}_n (\mu_0 1^{n-|\mu_0|})$, with 
$|\mu_0| \leq 14$ (and, of course, the smallest part of $\mu_0$ larger than one) ($135$ cases altogether), see the
output file:

{\tt http://www.math.rutgers.edu/\~{}zeilberg/tokhniot/oSn1} .

II. {\bf Explicit Expressions for} $ {\bf \psi^{(2)}_n (\mu_0 1^{n-|\mu_0|}) }$

The classical case is well-known, and easy.

$$
\psi^{(2)}_n (1^n)= \frac{1}{n+1} {{2n} \choose {n}} \quad (=C_n) \quad .
$$
We have:
$$
\psi^{(2)}_n (2 1^{n-2})= 
{\frac {9-5\,n+{n}^{2}}{ \left( 2\,n-1 \right)  \left( 2\,n - 3 \right)  \left( n+1 \right) }}
{{2n} \choose {n}}  \quad ,
$$
$$
\psi^{(2)}_n (3 1^{n-3})= 
{\frac {48-11\,n+{n}^{2}}{ 4 \left( 2\,n-1 \right)  \left( 2\,n -3 \right)  \left( n+1 \right) }}
{{2n} \choose {n}} \quad ,
$$
Note the {\bf remarkable} (proved!) identity:
$$
\psi^{(2)}_n (3 1^{n-3})= \frac{1}{2}\phi^{(2)}_{n+2} (32 1^{n-3}) \quad .
$$
It may be interesting to find a `natural' reason for this `coincidence'.

We also have:
$$
\psi^{(2)}_n (4 1^{n-4})= 
{\frac {2100-1354\,n+299\,{n}^{2}-26\,{n}^{3}+{n}^{4}}{ 4 \left( 2\,n-1 \right)  \left( 2\,n -3 \right)  \left( 2\,n-5 \right)  \left( 2\,n-7 \right) 
 \left( n+1 \right) }}
{{2n} \choose {n}} \quad ,
$$
$$
\psi^{(2)}_n (22 1^{n-4})= 
{\frac {525-316\,n+89\,{n}^{2}-14\,{n}^{3}+{n}^{4}}{ \left( 2\,n-1 \right)  \left( 2\,n -3\right)  \left( 2\,n-5 \right)  \left( 2\,n-7 \right)  \left( n+1
 \right) }}
{{2n} \choose {n}}  \quad ,
$$
$$
\psi^{(2)}_n (5 1^{n-5})= 
{\frac {10080-4342\,n+659\,{n}^{2}-38\,{n}^{3}+{n}^{4}}{ 16 \, \left( 2\,n-1 \right)  \left( 2\,n -3 \right)  \left( 2\,n-5 \right)  \left( 2\,n-7 \right) 
 \left( n+1 \right) }}
{{2n} \choose {n}}  \quad ,
$$
$$
\psi^{(2)}_n (32 1^{n-5})= 
{\frac {2520-1045\,n+194\,{n}^{2}-20\,{n}^{3}+{n}^{4}}{ 4\, \left( 2\,n-1 \right)  \left( 2\,n -3 \right)  \left( 2\,n-5 \right)  \left( 2\,n-7 \right) 
 \left( n+1 \right) }}
{{2n} \choose {n}} \quad .
$$

For {\it all} the (proved!) explicit expressions for $\psi^{(2)}_n (\mu_0 1^{n-|\mu_0|})$, with 
$|\mu_0| \leq 14$ (and, of course, the smallest part of $\mu_0$ larger than one) ($135$ cases altogether), see the
output file:

{\tt http://www.math.rutgers.edu/\~{}zeilberg/tokhniot/oSn2} .

{\bf Sums over shapes with more rows}

For three and more rows, the sums are no longer closed-form, but, as we mentioned above, they always
satisfy a linear recurrence (alias difference) equation with polynomial coefficients, see the output 
files linked to in the above-mentioned webpage of this article.

{\bf Conclusion}

            On page 155 of [GPK] it says: 

{\it ``The numbers in Pascal's triangle satisfy, practically
speaking, infinitely many identities, so it is not too surprising that we can find some surprising
relationships by looking closely.''}

The aim of this note was to indicate that a similar statement seems to hold for the character tables of the symmetric groups Sn.
Just as importantly, it was a {\it case-study} in using a computer algebra system to prove deep identities,
way beyond the ability of mere humans.

\bigskip

{\bf  Acknowledgment}

We thank Richard Stanley,  Dennis Stanton, and John Stembridge for useful information.
Special thanks go to Shalosh B. Ekhad, for its extensive, very reliable,  computations.

{\bf References}

[AZ] M. Apagodu and D. Zeilberger,
{\it Multi-Variable Zeilberger and Almkvist-Zeilberger algorithms and the sharpening of Wilf-Zeilberger theory},
Adv. Appl. Math. {\bf 37} (2006),(Special issue in honor of A. Regev), 139-152; \hfill\break
{\tt http://www.math.rutgers.edu/\~{}zeilberg/mamarim/mamarimhtml/multiZ.html} \quad .

[BR] A. Berele and A. Regev, {\it Asymptotics of Young tableaux in the $(k,l)$ hook}, in: ``{\it Groups, Algebras, and Applications}'' 
(C.P. Milies, ed.), Contemporary Mathematics {\bf 537} (2011),  71-84. See also: {\tt http://arxiv.org/abs/1007.3833} \quad .

[BFK] F. Bergeron, L. Favreau and D. Krob, 
{\it Conjectures on the enumeration of tableaux of bounded height}, Discrete Math {\bf 139} (1995), 463-468.

[EZ] S. B. Ekhad and A. Regev,
{\it
Refined asymptotics and explicit recurrences for the numbers of Young tableaux in the $(k,l)$ hook for $k+l \leq 5$ },
Personal Journal of Shalosh B. Ekhad and Doron Zeilberger; July 28, 2010; \hfill\break
{\tt http://www.math.rutgers.edu/\~{}zeilberg/mamarim/mamarimhtml/hooker.html} \quad .

[Ge] I.  Gessel, {\it Symmetric functions and P-recursiveness},
Journal of Combinatorial Theory Series A {\bf 53} (1990), 257-285; \quad 
{\tt http://people.brandeis.edu/\~{}gessel/homepage/papers/dfin.pdf} \quad .

[Go] D. Gouyou-Beauchamps, 
{\it Chemins sous-diagonaux et tableaux de Young}, in: ``Combinatoire \'enum\'erative  (Montr\'eal 1985, G. Labelle and P. Leroux, eds.)'', 
Lect. Notes Math. {\bf 1234},  112-125,  Springer, 1986.

[GKP] R.L. Graham, D.E. Knuth and O. Patashnik, {\it ``Concrete Mathematics}'', 1st. ed., Addison-Wesley, 1989.

[M] I. G. Macdonald, {\it ``Symmetric Functions and Hall Polynomials''}, 2nd ed., Clarendon Press, Oxford, 1995.

[Re] A. Regev, {\it Asymptotic values for degrees associated with strips of Young diagrams},  Advances in Mathematics {\bf 41} (1981), 115-136.

[Ro] G. de B. Robinson, {\it On the representations of $S_n$}, Amer. J. Math. {\bf 60} (1938), 745-760.

[Sc] C. E. Schensted, {\it Largest increasing and decreasing subsequences}, Canad. J. Math {\bf 13} (1961), 179-191.

[Sl] N.J.A. Sloane, {\it The On-Line Encyclopedia of Integer Sequences}; {\tt http://oeis.org} \quad .

[Wh] D. White, {\it A bijection proving orthogonality of the characters of Sn},
Adv. in Math. {\bf 50} (1983),  160 - 186. 

[Wi] The Wikipedia Foundation, {\it Hook Length Formula}; \hfill\break
{\tt http://en.wikipedia.org/wiki/Hook\_length\_formula} \quad .

[Z1] D. Zeilberger, {\it A holonomic systems approach to special functions identities}, 
J. Comput. Appl. Math. {\bf 32} (1990), 321 - 368; \hfill \break
{\tt http://www.math.rutgers.edu/\~{}zeilberg/mamarim/mamarimPDF/holonomic.pdf} \quad .

[Z2] D. Zeilberger, {\it The number of ways of walking in $x_1 \geq \dots \geq  x_k \geq 0$ for $n$ days, starting and ending at the origin,
     where at each day you may either stay in place or move one unit in any direction,
         equals the number of $n$-cell standard Young tableaux with $\leq 2k + 1$ rows},
Personal Journal of Shalosh B. Ekhad and Doron Zeilberger, Dec. 6, 2007; \hfill\break
{\tt http://www.math.rutgers.edu/\~{}zeilberg/mamarim/mamarimhtml/lazy.html} \quad .
\bigskip
\hrule
\smallskip
Alon Regev, Rockford, IL, USA  \quad .
\smallskip
Amitai Regev, Department of Pure Mathematics, Weizmann Institute of Science, Rehovot 76100, Israel ;
amitai dot regev at weizmann dot ac dot il \quad .
\smallskip
Doron Zeilberger, Department of Mathematics, Rutgers University (New Brunswick), Hill Center-Busch Campus, 110 Frelinghuysen
Rd., Piscataway, NJ 08854-8019, USA.  \hfill\break
zeilberg at math dot rutgers dot edu \quad .
\smallskip
\hrule
\medskip
July 13, 2015.

\end